\documentclass[12pt]{article}

\usepackage{amssymb,amsmath,amsthm}
\usepackage{amsfonts}
\usepackage{enumerate}

\usepackage{fancyhdr}	
\usepackage{mdframed}

\usepackage{tikz}
\usepackage{tikz-cd}

\usepackage{geometry}
\newtheoremstyle{mytheoremstyle} 
    {10pt}                    
    {10pt}                    
    {\normalfont}                   
    {}                           
    {\bfseries}                   
    {.}                          
    {0.3cm}                       
    {}  
\theoremstyle{mytheoremstyle}

\theoremstyle{plain}

\fancypagestyle{plain}{%
  \fancyhf{}%
}

\newtheorem{theorem}{Theorem}[section]

\newtheorem{lemma}[theorem]{Lemma}
\newtheorem{proposition}[theorem]{Proposition}

\newtheorem{corollary}[theorem]{Corollary}
\newtheorem{problem}[theorem]{Problem}

\begin{document}

\title{Generalized tournament matrices with the same principal minors}
\author{
	A. Boussaïri\thanks{Corresponding author: A. Boussaïri. Email: aboussairi@hotmail.com}, A. Chaïchaâ, B. Chergui and S. Lakhlifi
}

\maketitle

\begin{abstract}
A generalized tournament matrix $M$ is a nonnegative matrix that satisfies $M
  +M^{t}=J-I$, where $J$ is the all ones matrix and $I$ is the identity matrix. 
  In this paper, a characterization of generalized tournament matrices with the 
  same principal minors of orders $2$, $3$, and $4$ is given. In particular, it 
  is proven that the principal minors of orders $2$, $3$, and $4$ determine the 
  rest of the principal minors.
\end{abstract}



\section{Introduction}
  Let $M=(m_{ij})$ be an $n\times n$ matrix. With each nonempty subset $X
\subseteq \{1,\ldots,n\}$, we associate the \emph{principal submatrix} $M[X]$ 
of $M$ whose rows and columns are indexed by the elements of $X$. A
\emph{principal minor} of $M$ is the determinant of a principal submatrix of $M$. 
The \emph{order} of a minor is $k$ if it is the determinant of a $k\times k$ 
submatrix. In this paper, we address the following problem.
\begin{problem}\label{prob:1}
  What is the relationship between matrices with equal corresponding 
principal minors.
\end{problem}
  Clearly, if two matrices are diagonally similar then they have the same 
corresponding principal minors. Conversely, it follows from the main result 
of Engel and Schneider \cite{engel1980matrices} that  two symmetric matrices 
with no zeroes off the diagonal having the same principal minors of order $1$, 
$2$ and $3$ are necessarily diagonally similar. 
Hartfiel and Lowey \cite{hartfiel1984matrices} identified a special class of  
matrices in which two matrices  with equal corresponding principal minors of 
all orders are diagonally similar up to transposition. This result was 
improved in \cite{boussairi2015skew} for skew-symmetric matrices with no 
zeroes off the diagonal by considering only the equality of corresponding 
principal minors of order $2$ and $4$.

  Boussaïri and Chergui \cite{boussairi2016transformation} consider the class 
of  skew-symmetric matrices with entries from $\{-1,0,1\}$ and such that all 
off-diagonal entries of the first row are nonzero. They characterize the 
pairs of matrices of this class that have equal corresponding principal 
minors of order $2$ and $4$. This characterization involves a new  
transformation that generalizes diagonal similarity up to transposition.  

A \emph{tournament matrix} of order $n$ is the adjacency matrix of some 
tournament. In other word, it is an $n\times n$ $(0, 1)$-matrix $M$ which 
satisfies

\begin{equation}\label{eq:1} M + M^{t} = J_n - I_n,\end{equation}
where  $J_n$ denotes the all ones $n\times n$ matrix and $I_n$ denotes the $n 
\times n$ identity matrix. Boussaïri et al. \cite{boussairi2004c3} 
characterize the pairs of  tournaments having the same $3$-cycles. Clearly, 
two tournaments have the same 3-cycles if and only if their adjacency 
matrices have the same principal minors of order $3$. This implies a 
characterization of tournament matrices with the same principal minors of 
order $3$.

A \emph{generalized tournament matrix} $M = (m_{ij})$ is a nonnegative matrix that 
satisfies \eqref{eq:1}. By definition $m_{ij} = 1 - m_{ji}\in [0, 1]$ for all $i
\neq j\in\{1,\ldots,n\}$. Thus, we can interpret $m_{ij}$ as the a priori
probability that player $i$ defeats player $j$ in a round-robin tournament
\cite{moon1970generalized}. 

In this work, we characterize the pairs of generalized tournament matrices with 
the same principal minors of order at most $4$. We prove in particular that if 
two generalized tournament matrices have the same principal minors of orders 
at most $4$, then they have the same principal minors of all orders.

\section{Preliminaries and main result}
  Let $T$ be a tournament with vertex set $V$. A \emph{clan} of $T$ is a subset 
$X$ of $V$, such that  for all $a, b\in X$ and $x\in V\setminus X$, $(a, x)$ is 
an arc of $T$ if and only if $(b, x)$ is an arc of $T$. For a subset $Y$ of $V
$, we denote by ${\rm Inv}(T, Y)$ the tournament obtained by reversing all the arcs 
with both ends in $Y$. If $Y$ is a clan, we call this operation \emph{clan 
reversal}. It is easy to check that clan reversal preserves $3$-cycles. 
Conversely, Boussaïri et al. \cite{boussairi2004c3} proved that two 
tournaments on the same vertex set have the same $3$ cycles if and only if 
one is obtained from the other by a sequence of clan reversals.

Let $M = (m_{ij})$ be an $n\times n$ matrix. A \emph{clan} $X$ of $M$ is a 
subset of $[n]:=\{1,\ldots,n\}$ such that for all $i, j\in X$ and $k\in [n]\setminus
X$, $m_{ik}=m_{jk}$ and $m_{ki} = m_{kj}$. Denote by $M[X, [n]\setminus X]$ the
submatrix  of $M$ whose rows and columns are indexed by elements of $X$ and $[n]
\setminus X$ respectively. Clearly, $X$ is a clan of $M$ if and only if $M[X,
[n]\setminus X] = \mathbf{1}\cdot v^{t}$ and $M[[n]\setminus X, X] = w\cdot
\mathbf{1}^{t}$ for some column vectors $v$ and $w$.
The empty set, the singletons $\{i\}$ where $i\in [n]$, and $[n]$ are clans
called \emph{trivial}. We say that $M$ is \emph{indecomposable} if all its clans are 
trivial, otherwise it is called \emph{decomposable}. For a subset $Y$ of $[n]$, 
we denote by ${\rm Inv}(M, Y)$ the matrix obtained from $M$ by replacing the entry $
m_{ij}$ by $m_{ji}$ for all $i, j\in Y$. As for tournaments, if $Y$ is a clan 
of $M$, we call this operation clan reversal.

Let $M$ be a tournament matrix and let $T$ be its corresponding tournament. A
subset  $X$ of $[n]$ is a clan of $M$ if and only if it is a clan of $T$. Moreover, 
for every $Y\subset [n]$, the corresponding tournament of ${\rm Inv}(M, Y)$ is ${\rm Inv}(
T, Y)$. As the two possible tournaments on $3$ vertices have different
determinants, we can write Theorem 2 of \cite{boussairi2004c3}
as follows.

\begin{theorem}\label{theo11}
  Let $A$ and $B$ be two tournament matrices. The following assertions are 
equivalent:
  \begin{enumerate}[i)]
    \item $A$ and $B$ have the same principal minors of order $3$.
    \item There exists a sequence $A_0 = A, \ldots, A_m=B$, such that $A_{i+1}
 = {\rm Inv}(A_i, X_i)$ where $X_i$ is a clan of $A_i$ for all $i\in\{0, \ldots, m-1
\}$.
  \end{enumerate}
\end{theorem}

This theorem solves Problem \ref{prob:1} completely in the case of tournament 
matrices. Another result, in relation to our work, is the following theorem
due to Lowey \cite{loewy1986principal}.

\begin{theorem}\label{theo:lowey}
  Let $A, B$ be two $n\times n$ matrices. Suppose that $n\geq 4$, $A$ 
irreducible and for every partition of $[n]$ into two subsets $X, Y$ with $|X|
\geq 2$ and $|Y|\geq 2$, $rank(A[X, Y])\geq 2$ or $rank(A[Y, X])\geq 2$. If $A
$ and $B$ have equal corresponding minors of all orders, then they are 
diagonally similar up to transposition.
\end{theorem}

Let $M$ be an $n\times n$ generalized tournament matrix and let $X, Y$ be a
bipartition of $[n]$  with $|X|\geq 2$ and $|Y|\geq 2$. It is not hard to prove
that if $rank(M[X, Y]) \leq 1$ and $rank(M[Y, X])\leq 1$ then $X$ or $Y$ is a
nontrivial clan of $A$. It follows that an indecomposable generalized
tournament matrix satisfies the conditions of Theorem \ref{theo:lowey}.
Another fact is  that if two generalized tournament matrices are diagonally
similar, then they are equal. Then, from Theorem \ref{theo:lowey}, we have the
following proposition.

\begin{proposition}\label{propo:12}
  Let $A$ and $B$ be two $n\times n$ generalized tournament matrices. Suppose 
that $n\geq 4$ and $A$ is indecomposable. If $A$ and $B$ have equal 
corresponding minors of all orders, then $A=B$ or $A=B^{t}$.
\end{proposition}
  
  It follows from Theorem \ref{theo11} that it is enough to consider only 
principal minors of orders at most $3$ in the case of tournament matrices. 
This fact is not true for arbitrary generalized tournament matrices. Indeed, 
we will give in Section \ref{section:indec} two indecomposable $4\times 4$
matrices which have the same principal minors of orders $2$ and $3$, but do
not have the same determinant. Afterward, we will prove that Proposition
\ref{propo:12} still holds if we consider principal minors of orders at most
$4$. Then, we prove the following theorem.

\begin{theorem}\label{theo:1}
  Let $A$ and $B$ be two $n\times n$ generalized tournament matrices. The 
following assertions are equivalent:
  \begin{enumerate}[i)]
    \item $A$ and $B$ have the same minors of orders at most $4$.
    \item $A$ and $B$ have the same minors of every order.
    \item There exists a sequence $A_0=A,\ldots,A_m=B$ of $n\times n$ 
          generalized tournament matrices, such that for $k=0,\ldots,m-1$,
          $A_{k+1} = {\rm Inv}(A_k, X_k)$, where $X_k$ is a clan of $A_k$.
  \end{enumerate}
\end{theorem}

  It is worth noting that the proof of Theorem \ref{theo:lowey} in
\cite{loewy1986principal} uses tools from linear algebra. It seems hard to 
prove Theorem \ref{theo:1} in a similar fashion, even for $n=5$. We will use 
graph theoretic tools via a correspondence between generalized tournament matrices 
and weighted oriented graphs.

  Let $M=(m_{ij})$ be an $n\times n$ generalized tournament matrix. For all
$i\neq j \in [n]$, $m_{ij}$ is in $[0, 1]$, and $m_{ij} = m_{ji}$ if and only
if $m_{ij} = 1/2$. Then, we associate to $M$ a  weighted oriented graph 
$\Gamma_M$ with vertex set $[n]:=\{1,\ldots,n\}$, such that $(i, j)$ is an arc
with weight $m_{ij}$ if and only if $m_{ij} \in (1/2, 1]$. Conversely, let
$\Gamma$ be a weighted  oriented graph with vertex set $[n]$ and weights
in $(1/2, 1]$. We associate to $\Gamma$ a generalized tournament matrix $M=(m_
{ij})$, such that if $(i, j)$ is an arc then $m_{ij}$ is equal to the weight 
of $(i, j)$, and $m_{ij} = m_{ji} = 1/2$ if $(i, j)$ and $(j, i)$ are not 
arcs of $\Gamma$. 

  This correspondence between generalized tournament matrices and weighted 
oriented graphs allows us to use some techniques from \cite{boussairi2004c3}
in the proof of Theorem \ref{theo:1}.

\section{Decomposable and indecomposable weighted oriented graphs}
Let $\Gamma$ be a weighted oriented graph with vertex set $V$. We write $x 
\overset{\alpha}{\rightarrow} y$ if $(x, y)$ is an arc of $\Gamma$ with 
weight $\alpha$, and $x \cdots y$ if there is no arc between $x$ and $y$. 
Similarly, if $X$ and $Y$ are two disjoint subsets of $V$, we write $X \overset
{\alpha}{\rightarrow} Y$ if $(x, y)$ is an arc with weight $\alpha$ for every 
$x\in X$ and $y\in Y$. If $X={x}$ we simply write $x\overset{\alpha}{
\rightarrow} Y$ and $Y\overset{\alpha}{\rightarrow}x$ instead of $\{x
\}\overset{\alpha}{\rightarrow} Y$ and $Y\overset{\alpha}{\rightarrow}\{x\}$. 
The notations $X\cdots Y$, $x\cdots Y$ and $Y\cdots x$ are defined in the 
same way.

  A \emph{clan} of a weighted oriented graph $\Gamma$ with vertex set $V$ 
is a subset $X$ of $V$ such that for every $x\in V\setminus X$, either $x 
\cdots X$, $x\overset{\alpha}{\rightarrow} X$ or $X\overset{\alpha}{
\rightarrow}x$ for some weight $\alpha$. The empty set, the singletons $\{x\}$
where $x\in V$, and $V$ are clans called \emph{trivial}.  We say that $\Gamma$ is
\emph{indecomposable} if all its clans are trivial, otherwise it is called
\emph{decomposable}. The notion of clans was introduced, under different names,
for graphs, digraphs, and more generally $2$-structures
\cite{ehrenfeucht1990theory}. The next proposition gives some basic properties
of clans.

  \begin{proposition}\label{eq:clan_prop}
  Let $\Gamma$ be a weighted oriented graph with vertex set $V$. Let $X$, $Y
$ and $Z$ be subsets of $V$.
    \begin{enumerate}[i)]
      \item If $X$ is a clan of $\Gamma$, then $X\cap Z$ is a clan of
      $\Gamma  [Z]  $.

      \item If $X$ and $Y$ are clans of $\Gamma$, then $X\cap Y$ is a clan of
      $\Gamma$.

      \item If $X$ and $Y$ are clans of $\Gamma$,\ such that $X\cap
      Y\neq \emptyset$, then $X\cup Y$ is a clan of $\Gamma$.

      \item If $X$ and $Y$ are clans of $\Gamma$, such that $X\setminus
      Y\neq \emptyset$, then $Y\setminus X$ is a clan of $\Gamma$.

      \item If $X$ and $Y$ are clans of $\Gamma$, such that $X\cap Y=\emptyset
      $, then either  $X\overset{\alpha}{\rightarrow}Y$, $Y\overset{\alpha}{
      \rightarrow}X$ or $X\cdots Y$ for some weight $\alpha$.
    \end{enumerate}
  \end{proposition}
  
  The following theorem due to Ehrenfeucht and Rozenberg \cite{ehrenfeucht1990}
shows that indecomposability is hereditary.

  \begin{theorem}\label{theo:rozenberg}
    Let $\Gamma$ be an indecomposable weighted oriented graph with $n\geq5$ 
  vertices. Then, $\Gamma$ contains an indecomposable weighted oriented graph 
  with $n-1$ or $n-2$ vertices.
  \end{theorem}

  A weighted oriented graph $\Gamma$ is said to be \emph{separable} if its 
vertex set $V$ can be partitioned into two non empty clans, otherwise it is 
\emph{inseparable}. If $\Gamma$ is separable, then there exists a bipartition 
$X, Y$ of $V$ such that $X\overset{\alpha}{\rightarrow} Y$ for some weight $
\alpha$, or $X\cdots Y$. In the first case, $\Gamma$ is called
\emph{$\alpha$-separable}. We say that $\Gamma$ is \emph{$\alpha$-linear} if
its vertices can be ordered into a sequence $x_1,\ldots,x_n$ such that $x_i
\overset{\alpha}{\rightarrow} x_j$ if $i<j$. The notions defined above can be 
extended naturally to generalized tournament matrices.
  
  By definition, a tournament is inseparable if and only if it is irreducible.
It is well-known that every irreducible tournament with $n$ vertices contains
an irreducible tournament with $n-1$ vertices. The next theorem extends this
result to weighted oriented graphs and will be used in the proof of the main
theorem.
  
\begin{theorem}\label{theo:moon_wog}
  Let $\Gamma$ be an inseparable weighted oriented graph with $n\geq5$ 
vertices. Then, $\Gamma$ contains an inseparable weighted oriented graph 
with $n-1$ vertices.
\end{theorem}
  
\begin{proof}
  Suppose that $\Gamma$ is decomposable and let $C$ be a non trivial clan 
of $\Gamma$. Let $u$ be a vertex in $C$. We will prove that $\Gamma[V\setminus\{
u\}]$ is inseparable. Suppose, for the sake of contradiction, that $\Gamma[V
\setminus\{u\}]$ is separable, and let $X, Y$ be a bipartition of $V
\setminus\{u\}$ into two clans. Without loss of generality, we can suppose 
that $X\overset{\alpha}{\rightarrow} Y$. Since $C\neq V$, we have $C
\setminus\{u\}\neq X\cup Y$.
    \\
  \textbf{1.} If $C\setminus\{u\}\subseteq X$, then $C\setminus\{u\} \overset{
\alpha}{\rightarrow} Y$. As $C$ is a clan of $\Gamma$, $u 
\overset{\alpha}{\rightarrow} Y$. Hence, $X\cup\{u\} \overset{\alpha}{
\rightarrow} Y$, which contradicts the fact that $\Gamma$ is inseparable.
Similarly, $C\setminus\{u\}\subseteq Y$ yields a contradiction.\\
  \textbf{2.} If $C\setminus\{u\} \cap X$ and $C\setminus\{u\} \cap Y$ are 
non empty. Let $z\in V\setminus C$, we can suppose that $z\in X\setminus C$. 
We have $ z\overset{\alpha}{\rightarrow} Y$, in particular $ z\overset{\alpha
}{\rightarrow} Y\cap C$. Since $C$ is a clan, $ z\overset{\alpha}{\rightarrow}
 C$. It follows that $(X\setminus C)\overset{\alpha}{\rightarrow} V\setminus (
X\cap C)$.
    
  Suppose now that $\Gamma$ is indecomposable. The result is trivial if $
\Gamma$ contains an indecomposable graph with $n-1$ vertices. If no such 
graph exists, then, by Theorem \ref{theo:rozenberg}, there exist two distinct
vertices  $x, y\in V$ such that $\Gamma[V\setminus\{x, y\}]$ is indecomposable.
Then, $\Gamma[V\setminus\{x\}]$ or   $\Gamma[V\setminus\{y\}]$ is inseparable. 
Indeed, if $\Gamma[V\setminus\{x\}]$ is separable, then there exists a 
bipartition $X, Y$ of $V\setminus\{x\}$ into two clans. Suppose that $X
\overset{\alpha}{\rightarrow}Y$. If $V\setminus \{x, y\}\cap X$ and $V
\setminus \{x, y\}\cap Y$ are both non empty, then they are a 
bipartition of $V\setminus \{x, y\}$ into two clans, which contradicts the 
fact that $\Gamma[V\setminus \{x, y\}]$ is indecomposable. Hence, $V\setminus 
\{x, y\}$ is a clan of $\Gamma[V\setminus\{x\}]$. Similarly, if $\Gamma[V
\setminus\{y\}]$ is separable, then $V\setminus \{x, y\}$ is a clan of $\Gamma[
V\setminus\{y\}]$. It follows that if $\Gamma[V\setminus\{x\}]$ and $\Gamma[V
\setminus\{y\}]$ are both separable, then $V\setminus \{x, y\}$ is a clan of $
\Gamma$, which contradicts the assumption that $\Gamma$ is indecomposable.
\end{proof}

\section{Indecomposable generalized tournament matrices}\label{section:indec}
In this section, we improve Proposition \ref{propo:12} by showing that it is 
enough to consider principal minors of orders at most $4$. More precisely, we 
prove the following result. 
\begin{theorem}\label{propo:BILT}
   Let $A$ and $B$ be two $n\times n$ generalized tournament matrices. Suppose 
that $n\geq 4$ and $A$ is indecomposable. If $A$ and $B$ have equal 
corresponding principal minors of orders at most $4$, then $A=B$ or $A=B^{t}$.
\end{theorem}

  Let $A = (a_{ij})$ and $B = (b_{ij})$ be two $n \times n$ generalized 
tournament matrices with the same principal minors of orders $2$. Then, for all
$i\neq j\in [n]$, $a_{ij} = b_{ij}$ or $a_{ij} = 1- b_{ij}$. It follows that 
the set $\binom{ [n]}{2}$ can be partitioned into three subsets.
	\begin{itemize}
	\item $\mathcal{P}_{=}:=  \{    \{  i,j \}  \in \binom{ [  n]  }%
{2} : \text{ }a_{ij}=b_{ij}\text{ and }a_{ij}\neq1/2 \}  $
	\item $P_{\neq}:=  \{    \{  i,j \}  \in \binom{ [  n]  }{2} : \text{
}a_{ij}=1-b_{ij}\text{ and }a_{ij}\neq1/2 \}  $
  \item  $P_{1/2}:=  \{    \{  i,j \}  \in \binom{ [  n]  }{2} : \text{
}a_{ij}=b_{ij}=1/2  \}$
	
	\end{itemize}
	
  The \emph{equality graph} and the \emph{difference graph} of $A$ and $B$, 
denoted by $\mathcal{E}(A, B)$ and $\mathcal{D}(A, B)$ respectively, are the 
undirected graphs with vertex set $V=[n]$, and whose arc sets are $P_{=}$ and $
P_{\neq}$. It follows from the definition that \begin{align}\label{eq:transpose}
\mathcal{E}(A, B) = \mathcal{D}(A, B^{t}).\end{align}

In what follows, we give some information about generalized tournament 
matrices with the same principal minors of orders $2$ and $3$, via the 
equality and difference graphs.

\begin{lemma}\label{lemma:three_of_triangle}
  Let $A = (a_{i, j})$ and $B = (b_{i, j})$ be two $n \times n$ generalized 
tournament matrices with the same principal minors of orders $2$ and $3$.
For every $i,j,k\in  [  n]$ we have
  \begin{enumerate}[i)]

    \item if $  \{  i,j \}  \in P_{\neq}$ and
  $  \{  i,k \}  ,  \{  j,k \}  \in P_{=}$, then $a_{ik}=a_{jk}=b_{ik}=b_{jk}$.

    \item if $  \{
    i,j \}  \in P_{=}$ and $  \{  i,k \}  ,  \{  j,k \}  \in
    P_{\neq}$, then $a_{ik}=a_{jk}=1-b_{ik}=1-b_{jk}$.

    \item if $  \{
    i,j \}  \in P_{=}$ and $  \{  i,k \}  ,  \{  j,k \}  \notin
    P_{=}$, then $a_{ik}=1/2$  if and only if $a_{jk}=1/2$.

    \item if $  \{
    i,j \}  \in P_{\neq}$ and $  \{  i,k \}  ,  \{  j,k \}  \notin
    P_{\neq}$, then $a_{ik}=1/2$  if and only if $a_{jk}=1/2$.
  \end{enumerate}
\end{lemma}
	
		 \begin{proof}
    Let $i,j,k\in  [  n]  $. Then we have 
    \begin{align*}
      \det A[\{i, k\}] &= \det{B[\{i, k\}]}\\
      \det A[\{j, k\}] &= \det{B[\{j, k\}]}\\
      \det A[\{i, j, k\}] &= \det{B[\{i, j, k\}]}
    \end{align*}
    It follows that
    \begin{align}
       a_{i k} & = b_{ik} \mbox{ or } a_{ik} = 1 - b_{ik} \label{eq:01}\\
       a_{j k} & = b_{jk} \mbox{ or } a_{jk} = 1-b_{jk}\label{eq:02}\\
       a_{ik}-a_{ij}a_{ik}+a_{ij}a_{jk}-a_{ik}a_{jk} &=  b_{ik}-b_{ij}b_{ik}+b_{ij}b_{jk}-b_{ik}b_{jk} \label{eq:03}
    \end{align}
		
  If $  \{  i,j \}  \in P_{\neq}$ and $  \{  i,k \}  ,
  \{  j,k \}  \in P_{=}$, then $a_{ij}=1-b_{ij}$, $a_{ik}=b_{ik}$ 
and $a_{jk}=b_{jk}$. Using \eqref{eq:03}, we get $a_{jk}=a_{ik}$ and then $b_{
jk}=b_{ik}$. This proves assertion $i)$.
		
  To prove $iii)$ suppose that $  \{i,j \}  \in P_{=}$, $  \{  
j,k \}\notin P_{=}$ and $a_{ik}=1/2$. Then $b_{ij}=a_{ij} \neq 1/2$, $b_
{jk}=1-a_{jk}$ and $b_{ik}=1/2$.  By substituting in \eqref{eq:03}, we get $a_
{jk}=1/2$. Assertions $ii)$ and $iv)$  can be obtained from  $i)$ and $ii)$ 
by using  \eqref{eq:transpose}.
\end{proof}
	
\begin{proposition}\label{coro:same_clans}
  Let $A = (a_{i, j})$ and $B = (b_{i, j})$ be two $n \times n$ generalized 
tournament matrices. If $A$ and $B$ have the same principal minors of orders
$2$ and $3$, then the connected components of $\mathcal{E}(A, B)$ and
$\mathcal{D}(A, B)$ are clans of $A$ and $B$.
\end{proposition}

\begin{proof}
  By \eqref{eq:transpose}, it suffices to consider $\mathcal{E}(A, B)$. Let $C
$ be a connected component of $\mathcal{E}(A, B)$. If $C = [n]$ or $C = \{i\}$
 for some $i\in [n]$, then $C$ is a trivial clan of $A$ and $B$. Otherwise, 
let $i\neq j\in C$ be two adjacent vertices and let $k\in [n]\setminus C$. 
Then $  \{i,j \}  \in P_{=}$ and $  \{i,k \}  ,  \{j,k
 \}  \notin P_{=}$. We have to prove that $a_{ik} = a_{jk}$ and $b_{ik} 
= b_{jk}$. For this,  there are  two cases to consider.
  \begin{itemize}
    \item[1)] If $a_{ik} = 1/2$ then by assertion $iii)$ of Lemma
    \ref{lemma:three_of_triangle} we have $a_{jk} =  1/2$, and hence
    $b_{ik} = b_{jk} = 1/2$.
      
    \item[2)] If $a_{ik}\neq 1/2$, then $a_{jk}\neq 1/2$. It follows that $
      \{  i,k \}  ,  \{  j,k \}  \in P_{\neq}$. We conclude  
    by assertion $ii)$ of Lemma \ref{lemma:three_of_triangle}.
  \end{itemize}
\end{proof}

Let $A$ and $B$ be two $n \times n$ generalized tournament matrices with the 
same principal minors of orders $2$ and $3$. Suppose that $A$ is indecomposable.
By Proposition \ref{coro:same_clans}, $A=B$ or $A=B^{t}$ if 
and only if $\mathcal{E}(A, B)$ and $\mathcal{D}(A, B)$ are not both connected. 
In general, $\mathcal{E}(A, B)$ and $\mathcal{D}(A, B)$ can be both 
connected. Indeed, let $a,b\in  [0,1] \setminus \{1/2 \}$ and 
consider the matrix $M_{a, b}$ defined as follows
	
\[
  M_{a,b}=
  \begin{pmatrix}
  0 & a & b & b\\
  1-a & 0 & 1-a & b\\
  1-b & a & 0 & a\\
  1-b & 1-b & 1-a & 0
  \end{pmatrix}
\]
	
It easy to check that  the  matrices $M_{a,b}$ and $M_{1-a,b}$ have equal 
corresponding minors of orders $2$ and $3$, and that both $\mathcal{E}(M_{a,b}, M_{1
-a,b})$ and $\mathcal{D}(M_{a,b}, M_{1-a,b})$ are connected. Moreover, if $a
\neq b$ and $a\neq 1-b$, then $M_{a,b}$ and $M_{1-a,b}$ are indecomposable
and do not have the same determinant. This shows the necessity of the
equality of principal minors of order $4$ in the assumptions of Theorem
\ref{propo:BILT}. With this strengthening, we obtain the following result
which implies Theorem \ref{propo:BILT}.

\begin{proposition}\label{propo:several_classes}
Let $A$ and $B$ be two $n \times n$ generalized 
tournament matrices with the same principal minors of orders at most $4$.
If $A$ is inseparable, then  $\mathcal{E}(A, B)$ or $\mathcal{D}(A, B)$ is 
not connected.  
\end{proposition}

We will prove this proposition by induction on $n$. The next lemma allows
us to solve the base case $n=4$.
 	
\begin{lemma}\label{propo:connected_cases_gtm}
  Let $A = (a_{ij})$ and $B = (b_{ij})$ be two $4 \times 4$ generalized 
tournament matrices with equal corresponding minors of orders $2$ and $3$. If 
$\mathcal{D}(A, B)$ and $\mathcal{E}(A, B)$ are connected, then there exists 
a  permutation  matrix $P$ such that $A=PM_{a,b}P^{t}$ and  $B=PM_{1-a,b}P^{t
}$  where,  $a,b\in  [ 0,1] \setminus \{1/2 \}$. Moreover, $\det(A)=
\det(B)$ if and only if $a= b$ or $a= 1-b$.
\end{lemma}
	
\begin{proof}
 Suppose  that $\mathcal{D}(A, B)$ and $\mathcal{E}(A, B)$ are connected. The 
only possibility is that  $\mathcal{D}(A, B)$ and $\mathcal{E}(A, B)$ are 
disjoint paths of length  three. Then, there is a permutation matrix $P$ 
such that
\begin{itemize}
	\item the edges of $\mathcal{D}(P^{t}AP,P^{t}BP)$ are $\{1,2\}$,  $\{2,3\}$ and $\{3,4\}$.
	\item the edges of $\mathcal{E}(P^{t}AP,P^{t}BP)$ are $\{1,3\}$,  $\{1,4\}$ and $\{2,4\}$.
\end{itemize}

	Let $A^{'}:=P^{t}AP$ and $B^{'}:=P^{t}BP$. The off-diagonal entries of $A^{'
}=(a_{ij}^{'})$ and $B^{'}=(b_{ij}^{'})$ are not equal to $1/2$.  Moreover, 
we have $a_{12}^{\prime}=1-b_{12}^{\prime}$, $a_{32}^{\prime}=1-b_{32}^{\prime
}$, $a_{34}^{\prime}=1-b_{34}^{\prime}$, $a_{13}^{\prime}=b_{13}^{\prime}$, $a
_{14}^{\prime}=b_{14}^{\prime}$ and $a_{24}^{\prime}=b_{24}^{\prime}$. 
The matrices $A^{'}$ and $B^{'}$ have  equal corresponding minors of orders $2$ 
and $3$. Then, by assertions $i)$ and $ii)$ of Lemma \ref{lemma:three_of_triangle},
we get  $a_{12}^{\prime}=a_{32}^{\prime}=a_{34}^{\prime}$ and $a_{
13}^{\prime}=a_{14}^{\prime}=a_{24}^{\prime}$. Let $a:=a_{12}^{\prime}$ and $
b:=a_{13}^{\prime}$. We have $a,b\in  [ 0,1] \setminus \{1/2 \}$, $A
^{'}=M_{a,b}$ and $B^{'}=M_{1-a,b}$. Hence $A=PM_{a,b}P^{t}$ and  $B=PM_{1-a,b
}P^{t}$. Then $\det  (  A)  -\det  (  B)  =\det  (M_{
a,b})  -\det  (  M_{1-a,b})  = (  2b-1)   (2a-1
)   (  a-b)   (  a+b-1)$. It follows that $\det(A)=
\det(B)$ if and only if $a= b$ or $a= 1-b$.
\end{proof}

\begin{proof}[Proof of Proposition \ref{propo:several_classes}]
The result is trivial if $n=2$ or $n=3$. For $n=4$, suppose that $\mathcal{D}
(A,B)  $ and $\mathcal{E} (  A,B)  $ are both 
connected. By Lemma \ref{propo:connected_cases_gtm}, 
$A=PM_{a,a}P^{t}$ for some real number $a\in  [
0,1] \setminus \{1/2 \}$ and permutation matrix $P$.
Hence $A$ is separable. 

We will continue by induction on $n$ for $n\geq5$. Suppose, by contradiction, 
that $A$ is inseparable and that $\mathcal{D} (
A,B)$ and $\mathcal{E} (A,B)$ are connected. By Theorem
\ref{theo:moon_wog}, there is $i\in  [  n]$ such that the principal 
submatrix $A[ [n]  \setminus \{i\}]$ is inseparable. By induction 
hypothesis and without loss of generality, we can assume that $\mathcal{D}^{\prime}
:=\mathcal{D}(A[[n]\setminus i], B[[n]\setminus i])$ is not connected.
Since $\mathcal{E} (  A,B)
  $ is connected, there exists $j\neq i$ such that $  \{  i,j 
\}  $ is an edge of $\mathcal{E} (  A,B)  $. Let $C$ be the 
connected component of $\mathcal{D}^{\prime}$ containing $j$. As $A[ [  n
]\setminus \{i\}]$ is inseparable, there exists $k\in C$ and $h\in  
(   [n]  \setminus \{i\} )  \setminus C$ such that $a_{ij}\neq
 a_{hk}$. Let $C^{\prime}$ be the connected component of $\mathcal{D}^{\prime
}$ containing $h$. Since $\mathcal{D} (  A,B)  $ is 
connected, there exists $l\in C^{\prime}$ such that $  \{  l,i 
\}$ is an edge of $\mathcal{D} (  A,B)  $. By Proposition
\ref{coro:same_clans}, $C$ and $C^{\prime}$ are clans of $A[ [  n]\setminus\{i\}]$
and $B[ [  n]  \setminus \{i\}]$, then $a_{lj}
=a_{hk}$. It follows that $\det A[\{i,j,l\}]-\det B[\{i,j,l\}]= (  2a_{il}
-1)   (a_{ij}+a_{jl}-1)  = (  2a_{il}-1)   (  a_
{ij}-a_{lj})$, then $\det A[\{i,j,l\}]\neq \det B[\{i,j,l\}]$, because $
a_{il}\neq1/2$ and $a_{ij}\neq a_{lj}$. This contradicts the fact that $A$ 
and $B$ have the same principal minors of order $3$.
\end{proof}

  Let $A$ and $B$ be two generalized tournament matrices with the same principal
minors of orders at most $4$. Suppose that $A$ is indecomposable. Then $A$ is
inseparable. Moreover, by the proposition we have just proved, $\mathcal{E}(A, B)$
or $\mathcal{D}(A, B)$ is not connected. Since the connected components of
$\mathcal{E}(A, B)$ are $\mathcal{D}(A, B)$ are intervals of $A$, and the later
are all trivial, one of $\mathcal{E}(A, B)$ and $\mathcal{D}(A, B)$ is an empty
graph. Hence, $A=B$ or $A=B^{t}$. This proves Theorem \ref{propo:BILT}.

Let $\mathcal{F}$ be the family  of $4 \times 4$ matrices permutationally 
similar   to a matrix  $M_{a,b}$ with $a,b\in  [ 0,1] \setminus \{1/
2 \}$, $a\neq b$ and $a\neq 1-b$. We say  that a  matrix is \emph{$\mathcal{F
}$-free} if  it  contains no member  of $\mathcal{F}$ as a principal submatrix.
Let $A $ and $B $ be two $n \times n$  $\mathcal{F}$-free generalized 
tournament matrices with the same principal minors of orders at most $3$.
By Lemma \ref{propo:connected_cases_gtm}, $A$ and $B$ have the same principal
minors of orders at most $4$. Hence, for $\mathcal{F}$-free generalized
tournament matrices, it is enough to consider equality of principal minors of
orders at most $3$ in Theorem \ref{propo:BILT}.

\section{Proof of the main theorem}
  Let $A$ be a generalized tournament matrix and let $X$ be a clan of $A$. 
If $X$ is a trivial clan then ${\rm Inv}(A, X) = A$ or ${\rm Inv}(A, X) = A^{t}$.
In both cases $\det A = \det {\rm Inv}(A, X)$. Assume now that $X$ is a 
nontrivial clan of $A$. Then, up to permutation, $A$ can be written as follows
  \[
    A=
    \begin{pmatrix}
    A_{11} & \alpha\beta^{t}\\
    \beta\alpha^{t} & A[X]
    \end{pmatrix}
    \mbox{ and }
    {\rm Inv}(A, X) = 
    \begin{pmatrix}
    A_{11} & \alpha\beta^{t}\\
    \beta\alpha^{t} & A[X]^{t}
    \end{pmatrix}
    \mbox{,}
	\]
	
	where $\beta=
	\begin{pmatrix}
	1\\
	\vdots\\
	1
	\end{pmatrix}$ and 
  $\alpha=
	\begin{pmatrix}
	a_1\\
  a_2\\
	\vdots\\
	a_{n-|X|}
	\end{pmatrix}$, where $a_i\in [0, 1]$. By Proposition 3 of
  \cite{bankoussou2019spectral}, $A$ and ${\rm Inv}(A, X)$ have the same determinant.
  Then we have the following result.
   
  \begin{lemma}\label{lemma:12}
    Clan inversion preserves principal minors.
  \end{lemma}

It follows from this lemma that matrices obtained from a series of clan 
inversions have the same principal minors. This proves the implication $iii)
\Rightarrow ii)$ of Theorem \ref{theo:1}. The implication $ii)\Rightarrow iii)
$ is trivial. The remaining of the section is devoted to proving the implication
$i)\Rightarrow iii)$, that is, pairs of matrices with the same principal minors
of orders at most $4$ are obtained by a series  of clan inversions. 

We start by reducing the problem to the case of matrices with a common
nontrivial clan. 
\begin{proposition}\label{propo1}
  Let $A $ and $B $ be two $n \times n$ generalized tournament matrices with 
the same principal minors of orders at most $4$. Suppose that $A$ is 
inseparable. If $A$ and $B$ have no common nontrivial clans then $A = B$ or
$A = B^{t}$.
\end{proposition}

\begin{proof}
Since $A$ is inseparable, by Proposition \ref{propo:several_classes},
$\mathcal{E}(A, B)$ or $\mathcal{D}(A, B)$ is not connected. If $\mathcal{D}(A,
B)$ is not connected, then by Proposition \ref{coro:same_clans} its connected
components are common clans of $A$ and $B$. If $A$ and $B$ have no common
nontrivial clans, then the connected components of $\mathcal{D}(A, B)$ must be
singletons and hence $A = B$. If $\mathcal{E}(A, B)$ is not connected, using
\eqref{eq:transpose}, we get $A=B^{t}$.
\end{proof}

\begin{proposition}\label{propo22}
  Let $A$ and $B$ be two $n\times n$ generalized tournament matrices. If $A$ and
$B$ are $\alpha$-linear for some $\alpha>1/2$, then there exists a clan $X$ of 
$A$ such that ${\rm Inv}(A, X)$ and $B$ have a common nontrivial clan.
\end{proposition} 

\begin{proof}
  Let $\Gamma_A$ and $\Gamma_B$ be the corresponding graphs of $A$ and $B$. 
Without loss of generality, we can suppose that for all $i\neq j\in [n]$, $i
\overset{\alpha}{\rightarrow}j$ in $\Gamma_A$ if $i<j$. There 
exists a permutation $\sigma$ of $[n]$, such that for all $i\neq j\in [n]$,
$\sigma(i)\overset{\alpha}{\rightarrow}\sigma(j)$ in $\Gamma_B$ if $i<j$.
Consider the clan $X = \{1, \ldots, \sigma(1)\}$ of $A$. Clearly, $\sigma(1)
\overset{\alpha}{\rightarrow}[n]\setminus\{\sigma(1)\}$ in the graph
corresponding to ${\rm Inv}(A, X)$. Hence, $\{\sigma(2),\ldots,\sigma(n)\}$ is a
common nontrivial clan of ${\rm Inv}(A, X)$ and $B$. 
\end{proof}

  For the last case, when $A$ is separable and there is no $\alpha>1/2$ such 
that $A$ and $B$ are $\alpha$-linear, we need the following results.

\begin{lemma}\label{lemma3}
  Let $A$ be an $n\times n$ decomposable generalized tournament matrix and let
$I$ be a nontrivial clan of $A$. Let $x\in I$, then $A$ is
inseparable if and only if $A[([n]\setminus I) \cup \{x\}]$ is inseparable.
\end{lemma}

\begin{proof}
  Suppose that $V:=[n]$ can be partitioned into two clans $X$, $Y$ of $A$. If
$((V\setminus I) \cup \{x\})\cap X$ and $((V\setminus I) \cup \{x\})\cap Y$ are 
nonempty, then they are a bipartition of $(V\setminus I) \cup \{x\}$ into two 
clans of $A[(V\setminus I) \cup \{x\}]$. Otherwise, suppose for example that
$((V\setminus I) \cup \{x\})\cap X$ is empty, then $X\subset I\setminus \{x\}$.
Hence $\{x\}$, $V\setminus I$ is a bipartition of $(V\setminus I) \cup \{x\}$
into two clans. In both cases $A[(V\setminus I) \cup \{x\}]$ is separable.

  Conversely, let $X, Y$ be a bipartition of $(V\setminus I) \cup \{x\}$ into 
two clans of $A[(V\setminus I) \cup \{x\}]$ and assume for example that
$x\in X$. Then $X\cup I, Y$ is a bipartition of $V$ into two clans of $A$ and,
hence, $A$ is separable.
\end{proof}

\begin{proposition}
  Let $A$ and $B$ be two $n\times n$ generalized tournament matrices with the 
same principal minors of orders at most $4$. Then $A$ is inseparable if and 
only if B is inseparable.
\end{proposition}

\begin{proof}
  We proceed by induction on $n$. For $n=3$ the result is trivial. Suppose 
that $A$ is inseparable. If $B = A$ or $B = A^{t}$ then $B$ is inseparable. 
Otherwise, by Proposition \ref{propo1}, $A$ and $B$ have a common nontrivial 
clan $I$. Let $x\in I$, then by Lemma \ref{lemma3}, $A[(V\setminus I)\cup \{x
\}]$ is inseparable and so is $B[(V\setminus I)\cup \{x\}]$ by induction 
hypothesis. It follows by Proposition \ref{propo1} again that $B$ is 
inseparable.
\end{proof}

\begin{corollary}\label{remark1}
  Under the assumptions of the previous lemma, for $\alpha>1/2$, $A$ is $
\alpha$-separable if and only if $B$ is $\alpha$-separable.
\end{corollary}

Clearly, if a matrix $A$ is $\alpha$-linear, then for every clan $I$ of $A$, $
A[I]$ is $\alpha$-separable. Conversely, by induction on $n$, we obtain the
following result.
\begin{lemma}\label{lemma2}
  If there exists $\alpha>1/2$ such that for every clan $I$ of $A$, $A[I]$ is 
$\alpha$-separable, then $A$ is an $\alpha$-linear.
\end{lemma}

\begin{proposition}\label{propo:19}
  Let $A$ and $B$ be two $n\times n$ generalized tournament matrices with the 
same principal minors of orders at most $4$. Suppose that $A$ is 
separable. If there is no $\alpha>1/2$ such that $A$ and $B$ are $\alpha$-linear,
then $A$ and $B$  have a common nontrivial clan.
\end{proposition}

\begin{proof}
  Let $\Gamma_A$ and $\Gamma_B$ be the corresponding graphs of $A$ and $B$. 
Suppose that there exists a bipartition $X, Y$ of $[n]$ such that $X\cdots Y$ 
in $\Gamma_A$. Clearly, $X\cdots Y$ in $\Gamma_B$. As $n\geq 3$, $X$ or $Y$ 
is a common nontrivial clan of $A$ and $B$. Suppose now that $A$ is
$\alpha$-separable for some $\alpha>1/2$. Let $\mathcal{J}_A$ be the set of clans
$I$ of $A$ such that $A[I]$ is not $\alpha$-separable, $\mathcal{J}_B$ is
defined similarly. Assume that $A$ or $B$ is not $\alpha$-linear. Then, by Lemma
\ref{lemma2}, $\mathcal{J}_A \cup \mathcal{J}_B$ is not empty. Let $I$ be an
element of $\mathcal{J}_A \cup \mathcal{J}_B$ with maximum cardinality and assume,
for example, that $I\in\mathcal{J}_A$. Consider the smallest clan $\tilde{I}$ of $B$ 
containing $I$. Clearly, $B[\tilde{I}]$ is not $\alpha$-separable. Indeed, if 
$X, Y$ is a bipartition of $\tilde{I}$ such that $X\overset{\alpha}{\rightarrow}
Y$ in $\Gamma_B$, then $I\subset X$ or $I\subset Y$ because, by Corollary 
\ref{remark1}, $B[I]$ is not $\alpha$-separable. This contradicts the minimality
of $\tilde{I}$ because $X$ and $Y$ are both clans of $B$. Then, $\tilde{I}\in 
\mathcal{J}_B$ and, hence, $\tilde{I} = I$ by maximality of the cardinality 
of $I$. It follows that $I$ is a common nontrivial clan of $A$ and $B$.
\end{proof}
  
Now we are able to complete the proof of Theorem \ref{theo:1}.
The implications $\mathbf{iii) \Rightarrow ii)}$ and $\mathbf{ii)\Rightarrow
i)}$ are already proven. The proof of the implication $\mathbf{i) \Rightarrow
iii)}$ is similar to that of \cite[Theorem~2]{boussairi2004c3}, and it will be
added in order for the paper to be self-contained.

Let $A$ and $B$ be two $n\times n$ generalized tournament matrices with the same
principal minors of orders at most $4$. We want to prove that $B$ is obtained
from $A$ by a sequence of clan inversions. For this, we proceed by induction on
$n$. The result is trivial for $n=2$. Assume that $n\geq 4$.
There is nothing to prove if $A = B$ or $A=B^{t}$.
Otherwise, by Propositions \ref{propo1}, \ref{propo22} and
\ref{propo:19}, we can suppose that $A$ and $B$ have a common nontrivial clan $X$.
Let $x\in X$ and denote by $U$ the set $(V\setminus X)\cup \{x\}$. By induction 
hypothesis, there exist matrices $S_0=A[U],\ldots,S_l=B[U]$ such that $S_{k+1
}={\rm Inv}(S_k, Y_k)$, where $Y_k$ is a clan of $S_k$ for all $k\in \{0,\ldots,l-1\}$.
For each  $i\in\{0,\ldots,l-
1\}$, the subsets $\tilde{Y}_i$ of $V$ is defined from $Y_i$ as $\tilde{Y}_i =
 Y$ if $x\notin Y_i$ and $\tilde{Y}_i = Y_i\cup X$ if $x\in Y_i$. Now, the 
sequence $(\tilde{S}_i)$ is defined by $\tilde{S}_0 = A$ and for all $i \in \{0, 
\ldots, l-1\}$, $\tilde{S}_{i+1} = {\rm Inv}(\tilde{S}_i, \tilde{Y}_i)$. Clearly, $
\tilde{S}_m[U] = B[U]$, $\tilde{S}_m[X] = A[X]$ or $A[X]^{t}$, and since $A[X]
$ and $B[X]$ have also the same principal minors of orders at most $4$, $
\tilde{S}_m[X]$ and $B[X]$ have the same principal minors of orders at most $4$. 
By the induction hypothesis, there are matrices $R_0=\tilde{S}
_m[X],\ldots,R_p=B[X]$ such that $R_{i+1} = {\rm Inv}(R_i, Z_i)$, where $Z_i$ is a 
clan of $R_i$. By considering $\tilde{R}_0=A^{''}$ and for  all $i\in\{0,\ldots,p-1\}$,
$\tilde{R}_{i+1} = {\rm Inv}(\tilde{R}_i, Z_i)$, it is obtained that $\tilde{R}_p = B$.

\section{Remarks and Questions}
\textbf{1.} Let $T$ be an tournament with vertex set $V$. We can associate to $T$
the $3$-uniform hypergraph $\mathcal{H}_T$ with vertex set $V$ whose hyperedges
are the $3$-subsets of $V$ that induce $3$-cycles in $T$. We call this hypergraph
the \emph{$C3$-structure} of $T$. Clearly, not every $3$-uniform hypergraph arises
as the $C3$-structure of some tournament. Linial and Morgenstern 
\cite{linial2016number} asked if the $C3$-structure of tournaments can be
recognized in polynomial time. Some progress on this problem has been made in
\cite{boussairi20203}.

  As the determinant of a tournament on $3$ vertices is $1$ if it is a $3$-cycle
and $0$ otherwise, Linial and Morgenstern's problem can be stated matricially
as follows. Does there exist a polynomial time algorithm that decides if a
vector $P\in\{0, 1\}^{\binom{n}{3}}$ arises as the principal minors of order $3$
of a tournament matrix. This problem can be generalized naturally to generalized
tournament matrices.
  
\begin{problem}
  Is there a polynomial time algorithm that decides whether a collection
$(P_\alpha)_{\alpha\in 2^{[n]}, 2\leq|\alpha|\leq4}$ of real numbers arises
as the principal minors of orders $2,3$ and $4$ of a generalized tournament
matrix?
\end{problem}

\textbf{2.} Let $n\geq4$ be an integer. Denote by $GT_n$ the set of all
$n\times n$ generalized tournament matrices and by $PM_n$ the set of collections
$(P_\alpha)_{\alpha\subset [n], 2\leq|\alpha|\leq4}$ of real numbers that arise
as the principal minors of orders $2,3$ and $4$ of $n\times n$ generalized
tournament matrices. Let $\phi:GT_n\rightarrow PM_n$ the map which associates to
each generalized tournament matrix the collection of its principal minors of
orders $2,3$ and $4$. By Theorem \ref{theo:1}, the determinant of a generalized
tournament matrix is determined by the principal minors of orders at most $4$.
Hence, there exists a unique map $\psi:PM_n\rightarrow\mathbb{R}$, such that
$\psi o\phi(M) = \det(M)$, for every $n\times n$ generalized tournament matrix
$M$. That is, the following diagram is commutative.

\begin{center}
  \begin{tikzcd}
    GT_n \arrow[r, "\phi"] \arrow[dr, "\det"]
    & PM_n \arrow[d, "\psi"]\\
    & \mathbb{R}
  \end{tikzcd}
\end{center}

We can ask if the map $\psi$ can be found explicitly, that is if the determinant
of an $n\times n$ generalized tournament matrix can be expressed in terms of its
principal minors of orders at most $4$.

\bibliographystyle{plain}
\bibliography{bibpaper}

\begin{thebibliography}{10}

\bibitem{bankoussou2019spectral}
Edward Bankoussou-mabiala, Abderrahim Boussa{\"\i}ri, Abdelhak
  Cha{\"\i}cha{\^a}, Brahim Chergui, and Soufiane Lakhlifi.
\newblock On the spectral reconstruction problem for digraphs.
\newblock {\em arXiv preprint arXiv:1910.13914}, 2019.

\bibitem{boussairi2015skew}
Abderrahim Boussa{\"\i}ri and Brahim Chergui.
\newblock Skew-symmetric matrices and their principal minors.
\newblock {\em Linear Algebra and its Applications}, 485:47--57, 2015.

\bibitem{boussairi2016transformation}
Abderrahim Boussa{\"\i}ri and Brahim Chergui.
\newblock A transformation that preserves principal minors of skew-symmetric
  matrices.
\newblock {\em The Electronic Journal of Linear Algebra}, 32:131--137, 2017.

\bibitem{boussairi20203}
Abderrahim Boussa{\"\i}ri, Brahim Chergui, Pierre Ille, and Mohamed Zaidi.
\newblock 3-uniform hypergraphs: Decomposition and realization.
\newblock {\em Contributions to Discrete Mathematics}, 15(1):121--153, 2020.

\bibitem{boussairi2004c3}
Abderrahim Boussa{\"\i}ri, Pierre Ille, G{\'e}rard Lopez, and St{\'e}phan
  Thomass{\'e}.
\newblock The ${C}_3$-structure of the tournaments.
\newblock {\em Discrete mathematics}, 277(1-3):29--43, 2004.

\bibitem{ehrenfeucht1990}
Andrzej Ehrenfeucht and Grzegorz Rozenberg.
\newblock Primitivity is hereditary for 2-structures.
\newblock {\em Theoretical Computer Science}, 70(3):343--358, 1990.

\bibitem{ehrenfeucht1990theory}
Andrzej Ehrenfeucht and Grzegorz Rozenberg.
\newblock Theory of 2-structures, part i: Clans, basic subclasses, and
  morphisms.
\newblock {\em Theoretical Computer Science}, 70(3):277--303, 1990.

\bibitem{engel1980matrices}
Gernot~M Engel and Hans Schneider.
\newblock Matrices diagonally similar to a symmetric matrix.
\newblock {\em Linear Algebra and its Applications}, 29:131--138, 1980.

\bibitem{hartfiel1984matrices}
DJ~Hartfiel and R~Leowy.
\newblock On matrices having equal corresponding principal minors.
\newblock {\em Linear algebra and its applications}, 58:147--167, 1984.

\bibitem{linial2016number}
Nati Linial and Avraham Morgenstern.
\newblock On the number of 4-cycles in a tournament.
\newblock {\em Journal of Graph Theory}, 83(3):266--276, 2016.

\bibitem{loewy1986principal}
Raphael Loewy.
\newblock Principal minors and diagonal similarity of matrices.
\newblock {\em Linear algebra and its applications}, 78:23--64, 1986.

\bibitem{moon1970generalized}
John~W Moon and NJ~Pullman.
\newblock On generalized tournament matrices.
\newblock {\em SIAM Review}, 12(3):384--399, 1970.

\end{thebibliography}

\end{document}